\documentclass{amsart}

\usepackage{amssymb,amsmath,amsthm,graphics,amscd,amsfonts}

\theoremstyle{plain}
\newtheorem{theorem}{Theorem}[section]
\newtheorem{proposition}[theorem]{Proposition}
\newtheorem{corollary}[theorem]{Corollary}
\newtheorem{lemma}[theorem]{Lemma}

\newtheorem{definition}[theorem]{Definition}

\newtheorem{remark}[theorem]{Remark}

\newcommand{\bfb}{{\bf b}}

\newcommand{\bfC}{{\mathbb C}}

\newcommand{\bfR}{{\mathbb R}}
\newcommand{\bfZ}{{\mathbb Z}}

\newcommand{\bfQ}{{\mathbb Q}}

\newcommand{\bark}{{\overline k}}

\newcommand{\barpartial}{{\overline \partial}}

\newcommand{\barz}{{\overline z}}
\newcommand{\barv}{{\overline v}}

\newcommand{\calH}{\mathcal{H}}

\newcommand{\tildeX}{\widetilde {X}}

\newcommand{\tildeomega}{{\widetilde \omega}}
\newcommand{\tildeeta}{{\widetilde \eta}}

\newcommand{\mapright}[1]{\smash{\mathop{   \hbox to 0.7cm{\rightarrowfill}}
  \limits^{#1}}}

\def\p{\partial}

\def\p{\partial}

\def\la{\langle}
\def\ra{\rangle}

\def\ba{\begin{array}}
\def\ea{\end{array}}

\begin{document}

\title
{Uniqueness and examples of compact toric Sasaki-Einstein metrics}
\author{Koji Cho}
\address{Department of Mathematics, Kyushu University,
6-10-1, Hakozaki, Higashiku, Fukuoka-city, Fukuoka 812-8581 Japan}
\email{cho@math.kyushu-u.ac.jp}
\author{Akito Futaki}
\address{Department of Mathematics, Tokyo Institute of Technology, 2-12-1,
O-okayama, Meguro, Tokyo 152-8551, Japan}
\email{futaki@math.titech.ac.jp}
\author{Hajime Ono}
\address{Department of Mathematics, Tokyo Institute of Technology, 2-12-1,
O-okayama, Meguro, Tokyo 152-8551, Japan}
\email{ono@math.titech.ac.jp}

\date{November 13, 2006 }

\begin{abstract} 
In \cite{FOW} it was proved that, given a compact
toric Sasaki manifold with positive basic first Chern class and
trivial first Chern class of the contact bundle, one can find a
deformed Sasaki structure
on which a Sasaki-Einstein metric exists.
In the present paper 
we first prove the uniqueness of such Einstein metrics on compact toric Sasaki manifolds 
modulo the action of the identity component of the automorphism
group for the transverse holomorphic structure, and secondly
remark that the result of \cite{FOW} 
implies the existence of compatible Einstein
metrics on all compact Sasaki manifolds obtained from the toric diagrams
with any height, or equivalently on all compact toric Sasaki manifolds
whose cones have flat canonical bundle. We further show that there exists an infinite 
family of inequivalent toric Sasaki-Einstein metrics
on $S^5 \sharp k(S^2 \times S^3)$ for each positive integer $k$.
\end{abstract}


\maketitle

\section{Introduction}
In \cite{FOW} the existence of an Einstein metric is proved on a compact toric
Sasaki manifold with positive basic first Chern class and trivial first Chern class of the
contact bundle $D$; These two conditions will be denoted by $c_1^B > 0$ and 
$c_1(D) = 0$. The purposes of the present paper
is firstly to prove the uniqueness of Sasaki-Einstein metrics up to a connected Lie group
action and secondly to clarify the meaning of the assumptions $c_1^B > $ and
$c_1(D) = 0$ in relation with
toric diagrams.

A Sasaki manifold 
is a Riemannian manifold $(S, g)$ 
whose cone manifold 
$(C(S), \overline g)$ with $C(S) \cong  S\times \bfR^+$ and $\overline g = dr^2  + r^2g$ 
is K\"ahler where $r$ is the standard coordinate on $\bfR^+$.
From this definition $S$ is odd-dimensional and we put $\dim S = 2m + 1$, and thus
$\dim_{C(S)} = m+1$.
A Sasaki manifold 
$(S, g)$ is said to be toric if the K\"ahler cone manifold $C(S)$ is toric, namely 
$(m+1)$-dimensional torus $G$ acts on $(C(S), \overline g)$ effectively as
holomorphic isometries. Note that $C(S)$ does not
contain the apex. Then $S$ is a contact manifold
with the contact form
$$\eta = (i(\barpartial - \p) \log r)|_{r=1}$$
where $S$ is identified with the submanifold $\{r=1\} \subset C(S)$, 
and has the Reeb vector field $\xi$ 
with the defining properties
$$i(\xi)\eta= 1\quad \mathrm{and}\quad i(\xi) d\eta = 0$$
where $i(\xi)$ denotes the inner product. 
The Reeb field $\xi$ is a Killing vector field on $S$ and also lifts to
a Killing vector field on $C(S)$, and thus $\xi$ is contained in the Lie algebra $\mathfrak g$
of $G$ since $G$ already has the maximal dimension of possible torus actions on $C(S)$.

The Reeb vector field $\xi$ generates a $1$-dimensional foliation, called the Reeb
foliation, on $S$. Since $\xi$ naturally lifts to a holomorphic vector field on $C(S)$ in the form 
$\xi - iJ\xi$ with 
$\xi = J(r\frac{\p}{\partial r})$
the Reeb foliation shares common local leaf spaces with
the holomorphic flow generated by $\xi - iJ\xi$ on $C(S)$. Thus the local leaf spaces give the Reeb foliation
a transverse holomorphic structure.
The contact structure of $S$ determines
a K\"ahler structure on the transverse holomorphic structure, which we call the transverse  
K\"ahler structure. 

Recall that a smooth differential form $\alpha$ on $S$ is basic if
$$ i(\xi)\alpha = 0\quad \mathrm{and}\quad \mathcal L_{\xi}\alpha = 0$$
where $\mathcal L_{\xi}$ denotes the Lie derivative by $\xi$.
The basic forms are preserved by the exterior derivative $d$ which decomposes into
$d = \p_B + \barpartial_B$, and we can define basic cohomology groups and
basic Dolbealt cohomology groups. We also have the transverse Chern-Weil theory and
can define basic Chern classes for complex vector bundles with basic transition functions.
The Sasaki manifold is said to have positive basic first Chern class if the first Chern class of the 
normal bundle of the Reeb foliation is represented by a positive basic $(1,1)$-form; as mentioned 
above this condition is denoted by $c_1^B > 0$. This is a necessary condition for the existence
of Sasaki-Einstein metric or equivalently the existence of positive transverse K\"ahler-Einstein metric.
There is another necessary condition $c_1(D) = 0$ as a de Rham cohomology class where
$D = \mathrm{Ker}\, \eta$ is the toric bundle. 
Coversely if $c_1^B > 0$ and $c_1(D) = 0$ then $c_1^B = \tau [d\eta]$ for some positive
constant $\tau$. See Proposition 4.3 in \cite{FOW} for more details.
Given a Sasaki manifold $(S, g)$, we say that another Sasaki metric $g^{\prime}$ is compatible with the Sasaki structure of $(S, g)$ if $g$ and $g^{\prime}$ have the same Reeb vector field and thus define the same transverse 
holomorphic structure.\par

The automorphism group of the the transverse holomorphic structure is the group of all
biholomorphic automorphisms of $C(S)$ 
which commute with the holomorphic flow generated by
$ \xi - iJ\xi$. Such automorphisms descend to an action on $S$ preserving the transverse
holomorphic holomorphic structure of the Reeb foliation, see section 2 for more detail.
In this paper we first prove
the uniqueness theorem of compatible Sasaki-Einstein metrics modulo
connected group actions of automorphisms
for the transverse holomorphic structure. 

\begin{theorem}\label{Main1}  Let $(S, g)$ be a compact toric Sasaki manifold with $c_1^B > 0$
and $c_1(D) = 0$. 
Then the identity component of the automorphism
group for the transverse holomorphic structure acts transitively 
on the space of all Sasaki-Einstein metrics compatible with $g$. 
\end{theorem}

In order to make clear which Sasaki manifolds the result of \cite{FOW} applies to, we wish to
explain the conditions $c_1^B > 0$ and $c_1(D) = 0$. Since a three dimensional 
Einstein manifold of positive scalar curvature is finitely covered by the standard three sphere
we may restrict ourselves to the case when the dimension of $S$ is bigger than or equal to
five.

\begin{theorem}\label{Main2} Let $S$ be a compact toric Sasaki manifold with 
$\dim S \ge 5$. 
Then the following three conditions are equivalent.
\begin{enumerate}
\item[(a)] $c_1^B > 0$ and $c_1(D) = 0$.
\item[(b)] The Sasaki manifold $S$ is obtained from a toric diagram with height $\ell$ for some
positive integer $\ell$ defined by $\lambda_1, \cdots, \lambda_d \in \mathfrak g$ and
$\gamma \in \mathfrak g^{\ast}$ (c.f. Definition \ref{good} and \ref{TD1}) and the Reeb field
$\xi \in \mathfrak g$ satisfies
$$ \la \gamma, \xi \ra = -m-1\ \ \mathrm{and}\ \ \la y, \xi \ra > 0\ \ \mathrm{for\ all}\ 
y \in C$$
where $C = \{y \in \mathfrak g^{\ast}| \la y, \lambda_j \ra \ge 0,\ j=1,\cdots, d\}$.
\item[(c)]  For some positive integer $\ell$, the $\ell$-th power 
$K^{\otimes\ell}_{C(S)}$ of the canonical line bundle $K_{C(S)}$
 is trivial.

\end{enumerate}
\end{theorem}

\begin{remark}\label{Gorenstein} We denote by 
$\overline{C(S)}$ the the closure of the cone $C(S)$, that is $C(S)$ plus the apex,
and consider it as an affine toric variety. 
It is a known fact that the condition of toric diagram with height $\ell$ is equivalent to
the apex being a $\bfQ$-Gorenstein singularity, that is the $\ell$-th power 
${\mathcal K}^{\otimes\ell}_{\overline{C(S)}}$ of the canonical sheaf 
${\mathcal K}_{\overline{C(S)}}$ is invertible, see \cite{Altman}.
\end{remark}

In the literature there are toric Sasaki manifolds denoted by $Y^{p,q}$ (\cite{1}, \cite{3}), 
$L^{p,q,r}$ (\cite{2}, \cite{3}), $X^{p,q}$ (\cite{4}) and
$Z^{p,q}$ (\cite{OYasui06}, \cite{Zpq}) which are constructed from toric diagrams with height $1$. 
They are all of positive basic first
Chern class by Theorem \ref{Main2}, and thus admit a Sasaki-Einstein metric by \cite{FOW}.
Combining the existence result of \cite{FOW} with Theorem \ref{Main1}  and Theorem \ref{Main2} we
get the following corollary. 

\begin{corollary}\label{Main3} Given a toric diagram, there is a unique
Sasaki structure whose cone is the one obtained from the toric diagram by
Delzant construction and on which there exist compatible Einstein metrics.  Moreover
the identity component of the automorphism group of the transverse holomorphic structure acts transitively on the
set of all compatible Einstein metrics.
\end{corollary}

Thus the Sasaki-Einstein metrics constructed in \cite{FOW} on $Y^{p,q}$ coincide
with those which have been known in the literature \cite{1}, \cite{3}.

Using diagrams we show that compact connected toric Sasaki manifolds associated
with toric diagrams of height bigger than $1$ are not simply connected and that
the converse is not true by giving an example.
We will also show the following.

\begin{theorem}\label{Main4}  For each positive integer $k$ there exists an infinite family of
inequivalent toric Sasaki-Einstein metrics on 
the $k$-fold connected sum $S^5 \sharp k(S^2\times S^3)$ of $S^2\times S^3$
with $S^5$.
\end{theorem}

The existence of (possibly non-toric) Sasaki-Einstein metrics on $S^5 \sharp k(S^2\times S^3)$
has been known by the works of Boyer, Galicki, Nakamaye and Koll\'ar (\cite{BoGaNa02}, 
\cite{BoGaKo05}, \cite{Kol04}), and that
the existence of toric Sasaki-Einstein metrics for all odd $k$'s has been known by van Coevering
(\cite{Coev06}).
Hence our result is new in that we obtain toric constructions for all even $k$'s. 
Moreover most of our examples should be
irregular while the previous ones are all quasi-regular.

We are grateful to Charles Boyer for pointing out our careless statement of the results
without the condition $c_1(D) = 0$ in the first version of the paper.

\section{Uniqueness of compatible Sasaki-Einstein metrics}

In K\"ahler geometry a well known method of proving uniqueness of constant scalar
curvature metrics is to use  geodesics on the space of all K\"ahler metrics in a fixed
K\"ahler class (\cite{mabuchi87-2}, \cite{donaldson98}, \cite{xxchen00}). This
idea becomes substantially simpler when the K\"ahler manifold under consideration is 
toric because the geodesic becomes a line segment expressed by the symplectic
potentials (\cite{guan99}). To prove Theorem \ref{Main1} we wish to use the same idea, 
but have to consider geodesics both on the space of transverse K\"ahler metrics on $S$
and on the space of K\"ahler metrics on $C(S)$. We therefore give an outline of the idea
in the case of compact K\"ahler manifolds first and then explain how we modify it in the Sasakian case.

Let $V$ be a compact K\"ahler manifold and $\calH$ the space of K\"ahler potentials
in a fixed K\"ahler class $[\omega_0]$:
$$ \calH = \{\varphi \in C^{\infty}(V)\ |\ \omega_{\varphi} = \omega_0 + \sqrt{-1} 
\p \barpartial \varphi > 0\}.$$
The tangent space $T_{\varphi}\calH$ at $\varphi \in \calH$ is identified with the
set $C^{\infty}(V)$ of all real smooth functions via
$$ \frac{d}{ds}_{|s=0} (\varphi + s \psi) = \psi \in C^{\infty}(V).$$
We have a natural Riemannian metric on $\calH$
$$ (\psi_1, \psi_2) := \int_M \psi_1\psi_2\, \omega_{\varphi}^n$$
where $n = \dim_{\bfC}V$ and $\psi_1, \psi_2 \in T_{\varphi}\calH \cong C^{\infty}(V)$.
For a smooth path $\varphi = \{\varphi_t\ |\ a \le t \le b\}$ in $\mathcal H$, 
let $\psi = \{\psi_t\ |\ a \le t
\le b \}$ be a vector field along $\varphi$, considered as
$$ \psi_t = \frac{d}{ds}_{|s=0} (\varphi_t + s\psi_t) \in T_{\varphi_t}\calH, \qquad a \le t \le b.$$
Then the covariant derivative by Levi-Civita connection is expressed as
\begin{equation}\label{covariant}
 \frac{D}{\p t} \psi = \dot{\psi_t} - \mathrm{Re} (\barpartial \dot\varphi_t, \barpartial \psi_t)_{\omega_t} = \dot{\psi_t} - \frac 12 (d\dot\varphi_t,d\psi_t)_{\omega_t} 
 \end{equation}
where $\omega_t = \omega_{\varphi_t}$.
Thus the equation of geodesics is given by
\begin{equation}\label{geodesic}
\ddot \varphi_t - |\barpartial \dot\varphi_t|_{\omega_t}^2 = 0.
\end{equation}
The K-energy, or Mabuchi energy, is defined by
$$ \mu(\omega_{\varphi}) := M(\omega_{\varphi}, \omega_0) = - \int_0^1 \left(
\int_V (\sigma_{\omega_t} - \overline{\sigma})\dot\varphi\,\omega_t^n \right)dt$$
where $\varphi_t := t\varphi$, $0\le t \le 1$, $ \omega_t = \omega_{\varphi_t}$ and
$$ \overline{\sigma} = \frac {\int_V n c_1(V) \omega_0^{n-1}}{\int_V \omega_0^n}.$$
The fundamental facts are the following:
\begin{itemize}
\item $\omega_{\varphi}$ is a critical point of $\mu$ if and only if $\omega_{\varphi}$
is a K\"ahler metric of constant scalar curvature.
\item The Hessian of $\mu$ is positive semi-definite, so $\mu$ is a convex function.
\item We have for any smooth path  $\{\varphi_t\}_{0 \le t \le 1}$
$$\frac{d^2\mu(\omega_t)}{dt^2} = \int_V |\barpartial Y_t|^2_{\omega_t} \,\omega_t^n 
- \int_V (\ddot{\varphi_t }- |\barpartial \dot{\varphi_t}|_{\omega_t}^2)\omega_t^n
$$
where $Y_t = \omega_t^{-1}(\barpartial \dot\varphi_t)$.
In particular, for geodesics $\varphi_t$ we have
\begin{equation}\label{second}
 \frac{d^2}{dt^2} \mu(\omega_t) = \int_V |\barpartial Y_t|^2_{\omega_t} \,\omega_t^n \ge 0.
\end{equation}
\end{itemize}

Now we can prove the uniqueness of constant scalar curvature metrics 
modulo the action of the identity component of the group of biholomorphic automorphisms of $V$
provided we have
a geodesic joining two such metrics as follows.
Suppose that both $\omega_0$ and $\omega_1$ are K\"ahler forms with constant scalar
curvature and that we have a geodesic $\omega_t$, $0 \le t \le 1$, joining them. Then it
follows from the above facts that 
$$ \frac{d^2\mu(\omega_t)}{dt^2} \ge 0, \qquad \frac{d\mu(\omega_t)}{dt} _{|t=0} = 0,
\qquad \frac{d\mu(\omega_t)}{dt} _{|t=1} = 0. $$
These imply $\frac{d\mu(\omega_t)}{dt} = 0 $ for all $t \in [0,1]$. But 
the equation (\ref{second}) shows
that $Y_t$ is a holomorphic vector field and $\omega_t$ is the pull-back of $\omega_0$
by an automorphism of $V$. 

The problem is then whether we can find a geodesic in the space of K\"ahler metrics. The geodesic 
equation is reduced to a degenerate Monge-Amp\`ere equation (\cite{semmes}).
The existence of $C^{1,1}$-solution was proved by X.X. Chen \cite{xxchen00}, but
$C^{1,1}$ geodesics are not enough to prove uniqueness and Chen used 
$\varepsilon$-approximations
of the solutions.
In the toric case, however, geodesics are obtained as a line segment of symplectic potentials
as shown by Guan \cite{guan99}, which is explained next.

Let $V$ be a toric K\"ahler manifold. Then $V$ is a completion of $(\bfC^{\ast})^n$ with
coordinates $w^1, \cdots , w^n$. Put $w^j = e^{z^j}$ and $z^j = x^j + i\theta^j$. Let
$F(x)$ be the K\"ahler potential of a $T^n$-invariant K\"ahler metric so that
$$ \omega = ig_{j\bark} dz^j \wedge d\barz^k = 
\frac i4 \frac{\p^2F}{\p x^j \p x^k} dz^j \wedge d\barz^k.$$
The {\it symplectic potential} $G$ is the Legendre transform of F:
$$ G(y) = \sum_{j=1}^n x^j \frac{\p F}{\p x^j} - F$$
with $y_j = \frac{\p F}{\p x^j}$. There is a symmetrical relation
$$ x^j = \frac{\p G}{\p y_j}, \qquad F(x) = \sum_{j=1}^n y_j \frac{\p G}{\p y_j} - G.$$
Thus as matrices
$$ (\frac{\p^2 F}{\p x^i \p x^j}) = (\frac{\p y_i}{\p x^j}) = (\frac{\p x^i}{\p y_j})^{-1} 
= (\frac{\p^2 G}{\p y_i \p y_j})^{-1}.$$

If $\{\omega_t\}$ is a curve in the space of K\"ahler forms and $F_t$ is the corresponding
K\"ahler potential then we have the $t$-dependent coordinates
$$y_t = \frac{\partial F_t}{\partial x}$$
on the image of the moment map.
Conversely if we start from a curve $G_t$ of symplectic potential with $t$-independent 
coordinates $y$ we have $t$-dependent coordinates 
$$x_t =  \frac{\partial G_t}{\partial y}$$
on $\bfR^n$.
To understand the geodesic equation better it is convenient to consider $F_t(x_t)$ in terms of
$t$-dependent coordinates $x_t$ and $G_t(y_t)$ in terms of $y_t$ with the relations
\begin{equation}\label{FxGy}
 y_{tj} = \frac{\p F_t}{\p x^j_t}, \qquad G_t(y_t) = \sum_{j=1}^n x_t^j \,\frac{\p F_t}{\p x_t^j} - F_t.
\end{equation}
We suppress $t$ for the notational convenience. First of all
$$ \frac {\p G(y)}{\p t} = \dot{G}(y) + \sum_{j=1}^n \frac {\p G}{\p y_j} \dot y_j
=  \dot{G}(y) + \sum_{j=1}^n x^j \dot y_j,$$
and
\begin{eqnarray*}
 \frac{\p}{\p t}(\sum_{j=1}^n x^j \frac{\p F}{\p x^j} - F)
 &=& \sum_{j=1}^n (\dot{x}^j \,\frac{\p F}{\p x^j} + x^j \,\frac{\p \dot{F}}{\p x^j} 
 + x^j\frac{\p^2 F}{\p x^j\p x^k} \dot{x}^k - \frac{\p F}{\p x^j}\dot{x}^j) - \dot{F} \\
 &=&  \sum_{j=1}^n x^j \dot y^j - \dot F(x).
\end{eqnarray*}
Thus 
\begin{equation}\label{dotGF}
 \dot G(y) = - \dot F(x).
 \end{equation}
 
 Taking the derivative of (\ref{dotGF}) we get
 \begin{equation}\label{ddotGF}
 \ddot{F} + \sum_{j=1}^n\frac{\p \dot{F}}{\p x^j}\dot{x}^j 
 = - \ddot{G} + \sum_{j=1}^n \frac{\p \dot{G}}{\p y^j}\dot{y}^j.
 \end{equation}
 In what follows we omit obvious indices and sum notations.
 Taking the derivative of $y = \frac{\p F}{\p x}$ we have
 \begin{equation}\label{ydot}
 \dot{y} = \frac{\p \dot{F}}{\p x} + \frac{\p y}{\p x}\dot{x}.
 \end{equation}
 
\begin{proposition}Let $M$ be a toric K\"ahler manifold.
\begin{enumerate}
\item[(a)] Let $F_t(x_t)$ and $G_t(y_t)$ be K\"ahler and 
symplectic potentials of $t$-dependent coordinates $x_t$ and $y_t$
satisfying the relations (\ref{FxGy}). Then
the geodesic equations are insensitive to $t$-dependent coordinates
in that 
$$ (\ddot{F} - \frac 12 |d\dot{F}|^2_t)(x_t) = 0$$
if and only if 
$$ \ddot{G}(y_t) = 0.$$
\item[(b)] For any two K\"ahler potential there exists a unique geodesic joining
them. 
In the action-angle coordinates $y,\ \theta$ with fixed standard symplectic form
$\omega = \sum_{i=1}^n dy_i \wedge d\theta^i$ the geodesic can be expressed
as $tG_1(y) + (1-t)G_0(y)$.
\end{enumerate}
\end{proposition}

\begin{proof} (a)\ \  It follows from  (\ref{ddotGF}), (\ref{ydot}) and (\ref{dotGF}) that
\begin{eqnarray*}
(\ddot{F} - \frac 12 |d\dot{F}|^2_t)(x_t) &=&
\ddot{F} - \frac{\p x}{\p y} \frac{\p \dot{F}}{\p x}\frac{\p \dot{F}}{\p x}\\
&=& - \ddot{G} - \frac{\p \dot{G}}{\p y}\dot{y} - \frac{\p \dot{F}}{\p x}\dot{x} 
- \frac{\p x}{\p y} \frac{\p \dot{F}}{\p x}\frac{\p \dot{F}}{\p x}\\
&=& - \ddot{G} - \frac{\p \dot{G}}{\p y}\dot{y} - \frac{\p \dot{F}}{\p x}\dot{x} 
-(\frac{\p x}{\p y}\dot{y} - \dot{x})\frac{\p \dot{F}}{\p x}\\
&=& - \ddot{G} - \frac{\p \dot{G}}{\p y} \dot{y} - \frac{\p \dot{F}}{\p y} \dot{y}\\
&=& - \ddot{G}(y_t).
\end{eqnarray*}
(b)\ \ The existence of a geodesic joining two K\"ahler potentials can be shown as
follows. We first fix $t$-independent coordinates $x$ on $\bfR^n$ and
$y$ on the convex polytope. Let $F_0(x)$ and $F_1(x)$ be two K\"ahler
potentials, and $y_0$, $y_1$, $G_0$ and $G_1$ be defined by 
$$y_0 = \frac{\p F_0}{\p x},\ G_0 = xy_0 - F_0\ ;\quad y_1 = 
\frac{\p F_1}{\p x},\ G_1 = xy_1 - F_1.$$
Put $y_t = ty_1 + (1-t)y_0$. Then $y_t$ is the moment map of the K\"ahler potential
$tF_1 + (1-t)F_0$ and thus gives coordinates on the image of the moment map. 
Put $G_t(y) = tG_1(y_t) + (1-t)G_0(y_t)$. Then obviously $\ddot{G}(y_t) = 0$, so the 
corresponding Legendre transform
$$F_t (x_t) = y_t\frac{\p G_t}{\p y_t} - G_t(y_t)$$
satisfies the geodesic equation $(\ddot{F} - \frac 12 |d\dot{F}|^2_t))(x_t) = 0$. Thus we get a geodesic $F_t(x)$ joining
$F_0$ and $F_1$. Notice that we inserted $t$-independent coordinates $x$ in $F_t$ so that
$F_t(x)$ becomes a geodesic in the original sense. 
We could perturb $y_t$ and get a different $x_t$, but $F_t(x)$ does not change because
of the uniqueness of the geodesic proved by X.-X. Chen \cite{xxchen00}. 
One can argue as in \cite{guan99} to
show that $F_t$ defines metrics on the whole K\"ahler manifold.

Taking the Legendre transform of 
the geodesic $F_t$ obtained in this way one sees that in the action-angle coordinates
on the polytope (see \cite{Abreu}) $G_t$ is expressed as $tG_1(y) + (1-t)G_0(y)$.
\end{proof}

Now we consider the case of transverse K\"ahler structure of compact Sasaki manifolds
of positive basic Chern class. We begin with the study of the automorphisms of
transverse holomorphic structure.

\begin{proposition} Let $S$ be a compact Sasaki manifold. Then the Lie algebra
of the automorphism group of transverse holomorphic structure is the Lie algebra
of all Hamiltonian holomorphic vector fields in 
the sense of Definition 4.4 of \cite{FOW}.
\end{proposition}
\begin{proof} Since the Reeb foliation has transverse holomorphic structure we can choose
local transverse holomorphic coordinates $z^1, \cdots, z^m$. They are used as part of 
local holomorphic coordinates as well as local coordinates on $S$. 
A local holomorphic vector field of the form $X^i \frac{\p}{\p z^i}$ is considered as a
local vector field on $C(S)$ as well as one on $S$. We will denote by $\tildeX^{\prime}$
the former and by $X^{\prime}$ the latter. Note that, along $\{r=1\}$, $X^{\prime}$ is 
the tangential part to $S \cong \{r=1\}$ of $\tildeX^{\prime}$. 

If a vector field $\tildeX$ generates a one-parameter group of 
automorphisms of $C(S)$ which commutes with the holomorphic flow generated by
$\xi - iJ\xi$ then $[\tildeX, \xi - iJ\xi] = 0$. If we set $\tildeX = Y - iJY$ with $Y$ the
real part of $\tildeX$ then $[\xi, Y] = 0$. From this one sees $[Y, J\xi] = J[Y, \xi] = 0$.
These mean that the holomorphic flow descends to $S$ and local leaf spaces, and that
$\tildeX$ descends to a holomorphic
vector field on each local leaf space. This local vector field can be regarded as a 
local vector field $\tildeX^{\prime}$ on $C(S)$ as well as $X^{\prime}$ on $S$.
Recall from \cite{FOW}
that the contact form $\eta$ on $S$
lifts to $C(S)$ as
$$\eta = 2d^c\log r = i(\barpartial - \p) \log r$$
where we use the same letter $\eta$ by the abuse of notation. 
We then have $\eta(\tildeX^{\prime}) = \eta(X^{\prime})$. This is because if $p : C(S) =
S \times \bfR_+ \to S$ is the projection then $\eta$ on $C(S)$ is $p^{\ast}$ and
$X^{\prime} = p_{\ast}\tildeX^{\prime}$.

Then $\tildeX$ can be expressed as
\begin{equation}\label{HH1}
 \tildeX = \eta(\tildeX) (\xi - iJ\xi) + (\tildeX^{\prime} - 
 \eta(\tildeX^{\prime})(\xi - iJ\xi)). 
 \end{equation}
Note that the right hand side is an orthogonal splitting. Taking $\barpartial$ of both sides of
(\ref{HH1}) we get
\begin{equation}\label{HH2}
\barpartial \eta(\tildeX) = \barpartial \eta(\tildeX^{\prime}) = \barpartial_B \eta(X^{\prime}).
\end{equation}
Taking the tangential 
component to $S$ of $\tildeX$ we obtain
\begin{equation*}
X := \eta(\tildeX)\xi + X^{\prime} - \eta(X^{\prime})\xi.
 \end{equation*}
Then since $\eta(X) = \eta(\tildeX)$, $X$ may be written as
\begin{equation*}
X = \eta(X)\xi + X^{\prime} - \eta(X^{\prime})\xi.
 \end{equation*}
Since $\eta$ is of the form 
$$ \eta =  dt - i\p_B f + i\barpartial_B f$$
where $t$ is the leaf coordinate with $\xi t = 1$ and $f$ is the K\"ahler potenitial
for the transverse K\"ahler form $\frac12 d\eta$ we have
$$ d\eta = 2i\p_B\barpartial_B f.$$
Hence we get 
$$ i(X)d\eta = i(X^{\prime})d\eta = 
2i(X^{\prime})\barpartial_B \eta = - 2\barpartial_B (\eta(X^{\prime})) = - 2\barpartial_B (\eta (X)).$$
 Hence $X$ is a Hamiltonian holomorphic vector field in 
the sense of Definition 4.4 of \cite{FOW}. It is easy to see that the Lie algebra
 consisting of all $\tildeX$ is isomorphic to the Lie algebra consisting of all
 Hamiltonian holomorphic vector fields $X$.
\end{proof}

Recall that a basic function $\varphi$ is a smooth function on $S$ such that
$\xi \varphi = 0$ where $\xi$ is the Reeb field.  
The transverse K\"ahler form $\omega^T$ is given by
$$ \omega^T = \frac12 d\eta$$
where 
$$\eta = 2d^c\log r|_{\{r=1\}\cong S} = i(\barpartial - \p) \log r|_{\{r=1\}\cong S} $$
and the transverse K\"ahler deformation is given by
$ \omega^T  + i\p_B\barpartial_B\varphi$ for some basic function $\varphi$
where $\p_B$ and $\barpartial_B$ are basic $\p$ and $\barpartial$-operators.
The tangent space to a transverse K\"ahler metric is therefore
the set of all basic functions $\varphi$. We may define geodesics in the space of transverse
K\"ahler metrics by the equation
$$ \ddot{\varphi} - |\barpartial_B \dot{\varphi}|_t^2 = 0.$$
We can derive the similar conclusion that if one can always find a geodesic joining two
K\"ahler potentials one can show that
the identity component of the automorphism group of the transverse holomorphic structure acts transitively on the space of
transverse K\"ahler metrics of constant scalar curvature by using the principle stated in the
Appendix of \cite{FOW}. In fact the corresponding equation to (\ref{second}) shows that
the geodesic joining two transverse K\"ahler metric of constant scalar curvature is
tangent to the Hamiltonian function of a Hamiltonian holomorphic vector field. 
Of course since transverse K\"ahler-Einstein metrics
have constant scalar curvature these arguments give the uniqueness of transverse K\"ahler-Einstein metrics modulo the action of the identity component of the automorphism group of 
the transverse holomorphic structure.

Suppose now that the compact Sasaki manifold $S$ is toric so that the cone $C(S)$ is
a toric K\"ahler manifold. We may also define the covariant derivative and geodesic
equation by (\ref{covariant}) and (\ref{geodesic}). By the above arguments we can always
find a geodesic joining two K\"ahler potentials on $C(S)$. 
The K\"ahler form on $C(S)$
is given by
$$ \omega = \frac 12 d(r^2\eta) = \frac 12 idd^c r^2 = \frac12 i\p\barpartial r^2.$$
A function on $S$ can be lifted to
$C(S) = S \times \bfR^+  $ and we use the same notation for a function on $S$ and
its lift to $C(S)$. 
The transverse K\"ahler deformation is given using a basic function $\varphi$ by
$$ \tildeeta = \eta + 2d_B^c\varphi = 2d^c \log (r\exp\varphi),$$
and hence the K\"ahler form $\omega$ on $C(S)$ is deformed by
\begin{equation}\label{Kdeform}
\tildeomega = \frac12 d(\tilde r^2\tildeeta) = \frac12 i\p\barpartial (r^2\exp(2\varphi)).
\end{equation}
Let $\mathcal K$ be the space of all K\"ahler metrics on $C(S)$ of 
the form $i\p\barpartial H$ for some real smooth function $H$ on $C(S)$, and 
$\mathcal K_{\omega}$
be the submanifold consisting of K\"ahler metrics obtained by transverse K\"ahler deformations
of the form (\ref{Kdeform}). 

\begin{lemma}\label{totallygeodesic} $\mathcal K_{\omega}$ is a totally geodesic
submanifold in $\mathcal K$. 
\end{lemma}
\begin{proof} By (\ref{Kdeform}) a curve in $\mathcal K_{\omega}$ is of the form
$\frac12 r^2\exp(2\varphi_t)$ so that its tangent vector is
$r_t^2\dot\varphi_t$ where we put $r_t = r\exp\varphi_t$. Similarly a vector filed
along $\frac12 r^2\exp(2\varphi_t)$ is of the form $r_t^2\psi_t$ for a curve
$\psi_t$ of basic functions. The covariant derivative of $r_t^2\psi_t$ along
$\frac12 r^2\exp(2\varphi_t)$ is computed by
\begin{eqnarray}\label{covariantK}
\frac{D}{dt}(r_t^2\psi_t) &=& 2r_t^2\dot\varphi_t \psi_t + r_t^2\dot\psi_t
- \frac12 (d(r_t^2\dot\varphi_t),d(r_t^2\psi_t))\\
&=& r_t^2(\dot\psi_t - \frac12(d\dot\theta_t,d\psi_t))\nonumber\\
&=& r_t^2\frac{D}{dt}\psi_t\nonumber
\end{eqnarray}
where the covariant derivative in the last term is the one for the transverse
K\"ahler structure. This shows that the covariant derivative of a vector
field in the tangent spaces of $\mathcal K_{\omega}$ along a curve in
$\mathcal K_{\omega}$ is tangent to $\mathcal K_{\omega}$. Thus
$\mathcal K_{\omega}$ is a totally geodesic submanifold.
\end{proof}
\begin{proposition}\label{geodesicK} 
A curve $i\p\barpartial(\frac12 r^2\exp(2\varphi_t))$ in $\mathcal K_{\omega}$
is a geodesic if and only if $\omega^T + i\p_B\barpartial_B\varphi_t$ is a geodesic in the 
space of the transverse K\"ahler 
metrics.
Moreover, 
for any given two transverse K\"ahler metrics with the same Reeb field $\xi$ corresponding
to toric K\"ahler cone metrics there exists a
unique geodesic joining them.
\end{proposition}
\begin{proof} The first statement follows immediately from (\ref{covariantK}).
Let $\omega^T_0$ and $\omega^T_1$ be the
transverse K\"ahler metrics with the common Reeb field $\xi$ corresponding to toric K\"ahler metrics 
$\omega_0$ and $\omega_1$ on $C(S)$.
Let $G_0$ and $G_1$ be the corresponding symplectic potentials. We use the
action-angle coordinates $y_i,\ \theta^i$.
Since $G_0$ and $G_1$ have common Reeb field $\xi$, $g = G_1 - G_0$
satisfies
$$ (\sum_{j=1}^{m+1} y_j \frac{\p}{\p y_j})\frac{\p g}{\p y_i}= 0$$
by (2.40) in \cite{MSY2}. Thus the geodesic 
$tG_1(y) + (1-t)G_0(y) = G_0 + tg$ 
joining 
$G_0$ and $G_1$ has the same Reeb field $\xi$, and thus the corresponding
K\"ahler potentials define the same transverse holomorphic structure.
From the first statement it follows that the geodesic $F_t = \frac12 r_t^2 = \frac 12 r^2\exp(2\varphi_t)$
with $\varphi_t$ basic smooth functions descends to a geodesic in the space of transverse K\"ahler
metrics.
\end{proof}

\noindent
{\it Proof of Theorem (\ref{Main1})}: 
Let $\omega^T_0$ and $\omega^T_1$ be the
transverse K\"ahler metrics corresponding to two Sasaki-Einstein metrics.
As proved in \cite{NT} and \cite{BGS} Lichnerowicz-Matsushima theorem for compact
K\"ahler manifolds of constant scalar curvature extends to compact Sasaki manifolds
of constant transverse scalar curvature we may assume that both $\omega^T_0$ and $\omega^T_1$ 
are invariant under the maximal compact subgroup of the group of automorphisms of 
the transverse holomorphic structure. In particular we may assume that they are invariant under
the maximal torus $G$, and thus we only need to consider the toric Sasaki-Eisntein metrics.

In \cite{MSY2} and \cite{FOW} it is shown that the volume functional of Sasakian structures
depends only on the Reeb fields, that 
there is a unique critical 
Reeb field $\xi$ which minimizes the volume functional and that only for the critical point $\xi$
the obstruction to the existence of transverse K\"ahler-Einstein metric vanishes.
Thus $\omega^T_0$ and $\omega^T_1$ must have a common
Reeb field $\xi$.
They can be joined by a geodesic by Proposition \ref{geodesicK}. We then apply the standard method
known in K\"ahler geometry as explained above, and it follows from 
(\ref{second}) that the geodesic is tangent to the Hamiltonian function of a Hamiltonian 
holomorphic vector field. This completes the proof.

\section{Toric diagrams}

We begin with the definition of a {\it good rational polyhedral cone}.

\begin{definition}[c.f. \cite{Lerman}] \label{good} Let $\mathfrak g^{\ast}$ be the dual of the Lie algebra $\mathfrak g$
of the $(m+1)$ dimensional torus $G$. Let $\bfZ_{\mathfrak g}$ be the integral lattice of $\mathfrak g$, that is
the kernel of the exponential map $\exp : \mathfrak g \to G$. 
A subset $C \subset \mathfrak g^{\ast}$ is a rational polyhedral cone if there exists a
finite set of vectors $\lambda_i \in \bfZ_{\mathfrak g}$, $1 \le i \le d$, such that
$$ C = \{ y \in \mathfrak g^{\ast}\ |\ \la y, \lambda_i \ra \ge 0\ \mathrm{for\ }\ i = 1, \cdots, d\}.$$
We assume that the set $\lambda_i$ is minimal in that for any $j$
$$ C \ne \{ y \in \mathfrak g^{\ast}\ |\ \la y, \lambda_i \ra \ge 0\ \mathrm{for\ all}\ i\ne j\}$$
and that each $\lambda_i$ is primitive, i.e. $\lambda_i$ is not of the form $\lambda_i = a\mu$
for an integer $a \ge 2$ and $\mu \in \bfZ_{\mathfrak g}$. 
(Thus $d$ is
the number of codimension $1$ faces if $C$ has non-empty interior.)
Under these two assumptions
a rational polyhedral cone $C$ with nonempty interior is good if the following condition holds.
If
$$ \{ y \in C\ |\ \la y, \lambda_{i_j}\ra = 0\ \mathrm{for\ all}\ j = 1, \cdots, k \}$$
is a non-empty face of $C$ for some $\{i_1, \cdots, i_k\} \subset \{1, \cdots, d\}$, then 
$\lambda_{i_1}, \cdots, \lambda_{i_k}$ are linearly independent
over $\bfZ$ and 
\begin{equation}\label{goodcondition}
\{ \sum_{j=1}^k a_j \lambda_{i_j}\ |\ a_j \in \bfR\} \cap \bfZ_{\mathfrak g} = 
\{ \sum_{j=1}^k m_j \lambda_{i_j}\ |\ m_j \in \bfZ\}.
\end{equation}
\end{definition}

Let $M$ be a $2m + 1$-dimensional compact connected contact toric manifold with the 
contact form $\eta$.
Namely there is an effective action of the $(m+1)$-dimensional torus $G$ which preserves
$\eta$. 
Then the moment map $\mu : M \to \mathfrak g^{\ast}$ is defined by
$$ \la \mu(p),X \ra = (\eta(X_M))(p)$$
where $X_M$ denotes the the vector filed on $M$ induced by $X \in \mathfrak g$.
We assume $\dim M = 2m + 1 \ge 5$. 
It is well-known (\cite{Lerman}) that if the action of $G$ is not free then the image
of the moment map is a good rational polyhedral cone.

\begin{definition}\label{TD1} An $(m+1)$-dimensional toric diagram with height $\ell$ is a
collection of $\lambda_i \in \bfZ^{m+1} \cong 
\bfZ_{\mathfrak g}$ satisfying $(\ref{goodcondition})$ and 
$\gamma  \in \bfQ^{m+1} \cong
(\bfQ_{\mathfrak g})^{\ast}$
 such that 
\begin{enumerate}
\item[(1)] $\ell$ is a positive integer such that $\ell\gamma$ is a primitive element of
the integer lattice $\bfZ^{m+1} \cong \bfZ^{\ast}_{\mathfrak g}$.
\item[(2)] $\la \gamma, \lambda_i\ra = -1$. 
\end{enumerate}
We say that a good rational polyhedral
cone $C$ is associated with a toric diagram of height $\ell$ if there exists a rational vector $\gamma$ 
satisfying $(1)$ and $(2)$ above.
\end{definition}
The reason why we use the terminology ``height $\ell$'' is because of the following proposition.
\begin{proposition}\label{TD4}  Using a transformation by an element of $SL(m+1,\bfZ)$ we may assume that
$$ \gamma = \left(\begin{array}{r} -\frac1\ell \\ 0 \\ \vdots \\ 0 \end{array}\right)$$
and the first component of $\lambda_i$ is equal to $\ell$ for each $i$.
\end{proposition}
\begin{proof}  
By elementary group theory there is an element $A$ of $SL(m+1, \bfZ)$ which sends the primitive
vector 
$\ell\gamma$ in $\bfZ^{m+1}$ to ${}^t(-1,0,\cdots,0)$ where the left upper $t$ denotes the transpose.
Then 
$A\gamma = {}^t(-\frac1\ell, 0, \cdots, 0).$
By transforming $\mathfrak g$ by ${}^tA^{-1}$, the transpose of $A^{-1}$, we get 
$$ \la A\gamma, {}^tA^{-1}\lambda_i \ra = \la \gamma , \lambda_i \ra = -1.$$
This implies the first component of ${}^tA^{-1}\lambda_i$ is $\ell$.
\end{proof}

Before we give a proof of Theorem \ref{Main2} we outline the proof of the following fact 
(c.f. \cite{Lerman}, \cite{MSY1}).
\begin{proposition}\label{TD3}
For each pair of a good rational
polyhedral cone $C$ and an element $\xi \in C_0^{\ast}$ where
$$C^{\ast}_0 = \{\xi \in \mathfrak g\ |\ \la v, \xi\ra > 0\ \mathrm{for\ all}\ v \in C\}$$
there is a compact connected toric Sasaki manifold $S$ whose moment map image is 
equal to $C \backslash \{0\}$ and whose
Reeb vector field is generated by $\xi$.
\end{proposition}
\begin{proof}[Outline of the proof] The construction of a contact manifold from 
a good rational polyhedral cone is the so-called Delzant construction.
Let $e_1, \cdots, e_d$ 
be the canonical basis of $\bfR^d$. Of course they
generate the lattice $\bfZ^d$. Let $\beta_{\bfZ} : \bfZ^d \to \bfZ_{\mathfrak g}\cong \bfZ^{m+1}$ 
be the homomorphism defined by
$$ \beta_{\bfZ}(e_i) = \lambda_i,$$
and $\beta_{\bfR} : \bfR^d \to \mathfrak g \cong \bfR^{m+1}$ be the natural linear map 
induced by $\beta_{\bfZ}$.
Since $C$ has non-empty interior then $\beta_{\bfR}$ is surjective, i.e. there is a subset 
$\{i_1, \cdots, i_{m+1}\}$ such that $\lambda_{i_1}, \cdots, \lambda_{i_{m+1}}$ are
linearly independent over $\bfR$. Then $\beta_{\bfZ}$ and $\beta_{\bfR}$ naturally induce a homomorphism
$ \beta_T : T^d \to G \cong T^{m+1} $
of the tori. Let $K$ be the kernel of $\beta_T$. We write $[a] \in T^d$ for the image of $a \in \bfR^d$.
Then 
$$ K = \{[a]\ |\ \sum_{i=1}^d a_i \lambda_i \in \bfZ_{\mathfrak g} \}.$$
It is a compact abelian subgroup of $T^d$ and its Lie algebra is $\ker \beta_{\bfR}$.
Note that $K$ is not connected in general. 
Let us consider the standard action of $T^d$ on $\bfC^d$ with the K\"ahler form
$\frac i2 \sum_{i=1}^d dv^i \wedge d\barv^j$ by
$$ [a] \cdot (v^1,\cdots,v^d) = (e^{2\pi i a_1}v^1,\cdots, e^{2\pi i a_d}v^d).$$
Consider the action of $K$ on $\bfC^d$ obtained as the restriction of the $T^d$-action
and the moment map $\mu_K : \bfC^d \to \mathfrak k^{\ast}$.
The K\"ahler cone manifold $C(S)$ is obtained as the K\"ahler quotient
$$ C(S) = (\mu_K^{-1}(0)\backslash \{0\})/K.$$
See \cite{Lerman} for more detail.
The closure $\overline{C(S)}$ is obtained as
$$ \overline{C(S)} = \mu_K^{-1}(0)/K $$
which is realized also as a normal complex analytic space via the standard method
using fans in algebraic geometry (\cite{Oda}). 

A Sasaki manifold $S$ is obtained as 
$$ S = (\mu_K^{-1}(0) \cap S^{2d-1})/K$$
where $S^{2d-1}$ is the standard $(2d-1)$-sphere in $\bfC^d$. 
This Sasaki metric is often called the canonical Sasaki metric, and 
the symplectic potential on $C(S)$ and 
the Reeb field are respectively given by
\begin{eqnarray*}
G^{can} &=& \frac 12 \sum_{i=1}^d l_i(y) \log l_i(y),\\
 \xi^{can} &=& \sum_{i=1}^d \lambda_i
\end{eqnarray*}
where $l_i(y) = \la \lambda_i, y \ra$. For a general 
Reeb field $\xi \in \mathfrak g$ a symplectic potential $ G^{can}_{\xi}$ on $C(S)$ is given by
$$ G^{can}_{\xi}(y) = \frac12 \sum_{i=1}^d l_i(y) \log l_i(y) + \frac12 l_{\xi}(y) \log l_{\xi}(y) - \frac12 l_{\infty}(y) \log l_{\infty}(y)$$
where $l_{\xi}(y) = \la \xi, y\ra$ and $l_{\infty} = \la \xi^{can}, y \ra$,
see \cite{MSY1} for more detail.
The corresponding
K\"ahler potential $F^{can}_{\xi}$, computed by the Legendre transform, is given by
$$ F^{can}_{\xi} = \frac 12 l_{\xi}(y),$$
see (61) in \cite{FOW}. Since the K\"ahler potential is equal to $\frac12 r^2$ then $r^2 = l_{\xi}(y)$ and
the Sasakian structure is determined via the identification $S \cong \{l_{\xi}(y) = 1\}
\subset C(S)$.
\end{proof}

\noindent
{\it Proof of Theorem \ref{Main2}}\ :\ First we prove that 
(a) implies (b). Suppose $c_1^B > 0$ and $c_1(D) = 0$.
By our assumption $(m+1)$-dimensional torus $G$ acts on $S$ preserving the Sasakian
structure. By Proposition 4.3 in \cite{FOW} we can then 
choose a $G$-invariant transverse K\"ahler form $\omega^T$ such that 
$$ c_1^B = (2m+2)[\omega^T].$$
Let $\rho^T$ be the Ricci form of $\omega^T$. 
Note that $\omega^T$ and $\rho^T$ are defined on each
local leaf spaces of the Reeb foliation, but they can be lifted to $S$ to define global 2-forms
on $S$. There exists a basic $G$-invariant smooth function $h$ on $S$ such that
\begin{equation}\label{MA1}
\rho^T = (2m+2)\omega^T + i\p_B\barpartial_B h
\end{equation}
on $S$. By an elementary curvature computation in Sasakian geometry the equation (\ref{MA1}) is
equivalent to
\begin{equation}\label{MA2}
\rho = - i \p\barpartial \log \det (F_{ij}) =i\p\barpartial h
\end{equation}
on $C(S)$ where $h$ is pulled back to $C(S) \cong \bfR_+ \times S$ so that $h$ satisfies
\begin{equation}\label{basic}
 r\frac r{\p r} h = \xi h = 0
\end{equation}
and where $F$ is the K\"ahler potential on $C(S)$,
$(e^{x^0 + i\theta^0}, \cdots, e^{x^m + i\theta^m})$ is the coordinates of 
$G^{\bfC} \cong (\bfC^{\ast})^{m+1}$ and
$$ F_{ij} = \frac{\p^2F}{\p x^i \p x^j}.$$
Note that since $F$ is $G$-invariant it is independent of $\theta^i$'s.
Since any $G$-invariant pluriharmonic function on $C(S)$ is an affine function then
there exists a $\gamma \in \mathfrak g^{\ast}$ such that 
\begin{equation}\label{MA3}
\log \det(F_{ij}) = - 2\sum_{i=0}^m\gamma_ix^i - h
\end{equation}
by replacing $h + \mathrm{constant}$ by $h$.
Using the Legendre transform $G$ of $F$ we get
\begin{equation}\label{MA4}
\log \det(G_{ij}) = 2\sum_{i=0}^m\gamma_i G_i + h.
\end{equation}
Using Abreu-Guillemin arguments about the boundary behavior of $G$ it is
shown in \cite{MSY1} that
\begin{equation}\label{MA5}
\la \lambda_j , \gamma \ra = -1\quad \mathrm{for}\ j = 1, \cdots, d.
\end{equation}
Since the moment map image has non-empty interior there are $(m+1)$
vectors $\lambda_{j_1}, \cdots, \lambda_{j_{m+1}}$ linearly independent
over $\bfR$. Hence one can consider $\gamma$ as a solution to the linear
equations
$$ 
\la \lambda_{j_i} , \gamma \ra = -1\quad \mathrm{for}\ i = 1, \cdots, m+1
$$
and sees that $\gamma \in \bfQ^{m+1}_{\mathfrak g^{\ast}}$. 
Choosing a positive integer $\ell$ such that $\ell\gamma$ is a primitive element
of the integer lattice. Since $\eta(\xi) = 1$ and the moment map on $C(S)$ is 
given by $\frac12 r^2\eta$ we have
$\la y, \xi \ra > 0$ for all $y \in C$. It is also shown in \cite{MSY1} that
$ \la \gamma,\xi \ra = -m-1$.
This proves that (a) implies (b). 

Next we prove that (b) implies (c). Return to the Delzant construction in the
proof of Proposition \ref{TD3}. 
One sees that $\mu_K^{-1}(0)$ is given by
$$ \mu_K^{-1}(0) = \{v \in \bfC^d\ |\ \sum_{i=1}^d b_i |v_i|^2 = 0\ \mathrm{
for\ all}\ b \in \mathfrak k \subset \bfR^d \}.$$
By Proposition \ref{TD4} we may assume that the first component of $\lambda_i$
is $\ell$ for each $i$. Recall that 
$$ \sum_{i=1}^d a_i\lambda_i \in \bfZ^{m+1}$$
for all $[a] \in K$. Looking at the first component we get
$$ \ell(a_1 + \cdots + a_d) \in \bfZ$$
for all $[a] \in K$. Thus
$$ (e^{2\pi i(a_1 + \cdots + a_d)} dv_1 \wedge \cdots \wedge dv_d)^{\otimes \ell}
= (dv_1 \wedge \cdots \wedge dv_d)^{\otimes \ell}.$$
Let $\bfb_1, \cdots, \bfb_{d-m-1}$ be a basis of $\mathfrak k$, and put
$\bfb_i = (b_{i1}, \cdots, b_{id})$ and
$$ X_i = \sum_{j=1}^d b_{ij} \frac \p{\p v^j}.$$
Then
\begin{equation}\label{topform}
 (i(X_1)\cdots i(X_{d-m-1}) dv_1 \wedge \cdots \wedge dv_d)^{\otimes \ell}
 \end{equation}
descends to a nowhere zero section of $K_{C(S)}^{\otimes \ell}$. 
Hence $K_{C(S)}^{\otimes \ell}$ is a trivial line bundle.
 This proves that (d) implies (c). Note that this proof shows the section (\ref{topform}) extends
 to the apex of $\overline{C(S)}$, as a token of $\bfQ$-Gorenstein property (c.f. Remark 
 \ref{Gorenstein}).

We now prove that (c) implies (a).
Suppose we are given a $G$-invariant Sasakian structure with Reeb field $\xi$ and with trivial
line bundle $K_{C(S)}^{\otimes\ell}$. Thus we have a $G$-invariant K\"ahler metric $\omega$ and a
nowhere vanishing holomorphic section $\Omega_1$ of $K_{C(S)}^{\otimes\ell}$. 
Let $h_1$ be defined by
$$ h_1 = \frac1\ell\log ||\Omega_1||^2$$
where the norm of $\Omega_1$ is taken with respect to $\omega$. 
Then the Ricci form $\rho$ of $\omega$ is written as
\begin{equation}\label{Ricci1}
\rho = \frac 1{2\pi} \p\barpartial h_1
\end{equation}
Let $h$ be the
average of $h_1$ by the action of $G$. Since $\rho$ is $G$-invariant we see from (\ref{Ricci1}) that
\begin{equation}\label{Ricci2}
\rho = \frac i{2\pi} \p\barpartial h.
\end{equation}
Starting from (\ref{MA2}) which is
identical to (\ref{Ricci2}) we get
(\ref{MA3}) and (\ref{MA4}) (though we do not have (\ref{basic})). Then 
it is shown in \cite{MSY1} that 
\begin{equation}\label{MA6}
\la \xi , \gamma \ra = -(m+1).
\end{equation}
The equation (\ref{MA3}) says that $e^h\det(F_{ij})$ is a flat metric on $C(S)$.
Consider the $(m+1)$-form $\Omega$ written as 
\begin{eqnarray}\label{MA7}
\Omega &=& e^{-i\sum_{i=0}^m \gamma_i \theta^i} e^{\frac h2} (\det(F_{ij}))^{\frac12}
dz^0 \wedge \cdots \wedge dz^m\\
&=& e^{-\sum_{i=0}^m \gamma_i z^i }dz^0 \wedge \cdots \wedge dz^m \nonumber
\end{eqnarray}
where we used  (\ref{MA3}). Then $\Omega$ is multi-valued if $\gamma$ is not integral
but only rational. Let $\ell_1$ be the positive integer such that $\ell_1\gamma$ is a primitive
element of the integer lattice. Then $\Omega^{\otimes\ell_1}$ is a holomorphic section 
of $K^{\otimes\ell_1}_{C(S}$ over the open set corresponding to the interior of the moment map
image. But since $||\Omega^{\otimes\ell_1}|| = 1$ we see that $\Omega^{\otimes\ell_1}$ extends
to the whole $C(S)$. We further have 
\begin{equation}\label{MA8}
\mathcal L_{\xi}\Omega = (m+1)i\Omega
\end{equation}
and 
\begin{equation}\label{e112}
\left(\frac{i}{2}\right)^{m+1}(-1)^{m(m+1)/2}\Omega\wedge \overline{\Omega}=
\exp (h)\frac{1}{(m+1)!}\omega^{m+1}.
\end{equation}
Since $\xi$ is decomposed into the holomorphic and the anti-holomorphic parts
$$\xi = \frac12(\xi - iJ\xi) + \frac 12 (\xi + iJ\xi) $$
with $J\xi = - r\frac \p{\p r}$ we have
\begin{eqnarray*}
\mathcal L_{\xi}\Omega &=& \mathcal L_{\frac12(\xi + ir\frac \p{\p r})}\Omega\\
&=& \frac{(m+1)i}2 \Omega + \frac i2 \mathcal L_{r\frac \p{\p r}}\Omega.
\end{eqnarray*}
From this and (\ref{MA8}) it follows that
$$ \mathcal L_{r\frac \p{\p r}}\Omega = (m+1)\Omega$$
and
$$ \mathcal L_{r\frac \p{\p r}}(\Omega\wedge\overline{\Omega}) = 2(m+1)\Omega\wedge
\overline{\Omega}.$$
On the other hand since $\omega = i\p\barpartial (r^2/2)$ we have
$$ \mathcal L_{r\frac \p{\p r}} \omega^{m+1} = 2(m+1)\omega^{m+1}.$$
Taking the Lie derivative of both sides of (\ref{e112}) by $r\frac \p{\p r}$ we get
$$ r\frac {\p h}{\p r} = 0. $$
Since $h$ is $G$-invariant we also have $\xi h = 0$. Hence (\ref{e112}) implies 
$[\rho^T] = (2m+2) [\omega^T]$ as basic cohomology classes, from which we get
$c_1^B > 0$ and $c_1(D) = 0$ by Proposition 4.3 in \cite{FOW}.
This completes the proof of Theorem \ref{Main2}.

\section{Examples and remarks on the fundamental groups}

As is mentioned in the introduction there are examples of 3 dimensional toric diagrams
of height 1 denoted by $X^{p,q}$, $Y^{p,q}$, $Z^{p,q}$ and $L^{p,q,r}$ 
known in physics literature and all the
corresponding Sasaki manifolds have Sasaki-Einstein metrics by 
the existence result of \cite{FOW}.
To check that these toric diagrams satisfy the goodness condition of Definition
\ref{TD1} the following proposition is useful.

\begin{proposition}\label{gooddiagram}   Let $C$ be a convex polyhedral cone in $\bfR^3$ given by
$$ C = \{ y \in \bfR^3\ |\ \la y, \lambda_i\ra \ge 0,\ j = 1, \cdots, d\}$$
with 
$$ \lambda_1 = \left(\begin{array}{c} 1 \\ p_1 \\ q_1\end{array}\right), \cdots,
\lambda_d = \left(\begin{array}{c} 1 \\ p_d \\ q_d\end{array}\right).$$
Then $C$ is good in the sense of Definition \ref{TD1} if and only if
either

(i)\ $|p_{i+1}- p_i | = 1 $ or $|q_{i+1} - q_i| = 1$ 

\noindent
or

(ii)\ $p_{i+1}- p_i $ and $q_{i+1} - q_i$ are relatively prime non-zero integers 

\noindent
for $i = 1,\ \cdots, d$ where we have put $\lambda_{d+1} = \lambda_1$.

Further the area of the 2-dimensional convex polytope formed by 
$$ \left(\begin{array}{c}  p_1 \\ q_1\end{array}\right), \cdots, 
\left(\begin{array}{c}  p_d \\ q_d\end{array}\right),  \left(\begin{array}{c}  
p_1 \\ q_1\end{array}\right)$$
is an invariant of the equivalent classes given by the action of some element of 
$SL(3,\bfZ)$ on the set of all
such diagrams.
\end{proposition}
\begin{proof}
Let $a_1$ and $a_2$ be real numbers such that $a_1 \lambda_i + a_2 \lambda_{i+1} \in \bfZ^3$.
Then we have 
\begin{equation}\label{prime1}
 a_1 + a_2 \in \bfZ,\ p_ia_1 + p_{i+1}a_2 \in \bfZ,\ q_i a_1 + q_{i+1}a_2 \in \bfZ.
 \end{equation}
It follows from these that
\begin{equation}\label{prime2}
 (p_{i+1} - p_i)a_2 \in \bfZ,\ (q_{i+1} - q_i)a_2 \in \bfZ.
 \end{equation}
If $(p_{i+1} - p_i)$ and $ (q_{i+1} - q_i)$ satisfy (i) or (ii) then there exist $s,\ t \in \bfZ$ 
such that $s(p_{i+1} - p_i)+ t(q_{i+1} - q_i) = 1$.  Then from (\ref{prime2}) we get 
$a_2 \in \bfZ$. From (\ref{prime1}) we also have $a_1 \in \bfZ$. Conversely if
$a_2$ satisfying (\ref{prime2}) is always in $\bfZ$ then $(p_{i+1} - p_i)$ and $ (q_{i+1} - q_i)$ 
are relatively prime. 

If a diagram of the first variable 1 is transformed to another by an element of $SL(3,\bfZ)$ then
the volume of the 3-dimensional truncated cone 
$$ \{ a_1\lambda_1 + \cdots + a_d\lambda_d \ |\ 0 \le a_1 \le 1, \cdots, 0 \le a_d \le 1, 0 \le a_1 + \cdots + a_d \le 1\}$$
is invariant. But this is equal to one thirds of the area described in the statement of the proposition.
This completes the proof of Proposition \ref{gooddiagram}.
\end{proof}

We give a simplest toric diagram of hight $\ell$. Let $C$ be the convex polyhedral cone defined by
$$ C = \{y \in \bfR^3\ |\ \la y, \lambda_i\ra \ge 0,\ j = 1,\ 2,\ 3\}$$
with
$$ \lambda_1 = \left(\begin{array}{c} 1 \\ 0 \\ 0\end{array}\right), 
\lambda_2 = \left(\begin{array}{c} 0 \\ 1 \\ 0\end{array}\right), 
\lambda_3 = \left(\begin{array}{c} 1 \\ 1 \\ \ell \end{array}\right).$$
Then this is a good cone and defines a smooth Sasaki manifold. One can show that
$$ \gamma =  \left(\begin{array}{c} -1 \\ -1 \\ \frac1\ell \end{array}\right).$$
Taking 
$$ A = \left(\begin{array} {rrr} 0 & 0 & -1 \\ -1 & 1 & 0 \\ 1 & 0 & \ell \end{array}\right)$$
we have 
$$ A\gamma = \left(\begin{array}{c} -\frac1\ell \\ 0 \\ 0\end{array}\right),\ 
{}^tA^{-1}\lambda_1 = \left(\begin{array}{c} \ell \\ 0 \\ 1\end{array}\right),\ 
{}^tA^{-1}\lambda_2 = \left(\begin{array}{c} \ell \\ 1 \\ 1\end{array}\right),\ 
{}^tA^{-1}\lambda_3 = \left(\begin{array}{c} \ell \\ 1\\ 2\end{array}\right).
$$
By following the Delzant construction one sees that the resulting Sasaki manifold
is the Lens space $S^5/\bfZ_\ell$. 

Next let $C$ be the convex polyhedral cone defined by
$$ C = \{y \in \bfR^3\ |\ \la y, \lambda_i\ra \ge 0,\ j = 1,\ 2,\ 3,\ 4\}$$
with
$$ \lambda_1 = \left(\begin{array}{c} 1 \\ 0 \\ 0\end{array}\right), 
\lambda_2 = \left(\begin{array}{c} 0 \\ 1 \\ 0\end{array}\right), 
\lambda_3 = \left(\begin{array}{c} 1 \\ 1 \\ \ell \end{array}\right),
\lambda_4 = \left(\begin{array}{c} 1 \\ 1 \\ \ell -1 \end{array}\right).
$$
Then this is a good cone and defines a smooth Sasaki manifold. The resulting
Sasaki manifold does not satisfy the conditions of Theorem \ref{Main2} because there
is no $\gamma$ with $\la \gamma, \lambda_j \ra = -1$ for $j= 1,\ 2,\ 3,\ 4$.

One can show that if $\ell > 1$ then the resulting Sasaki manifold $S$ is not simply
connected. This follows from a result of Lerman \cite{Lerman02} which is stated
as follows. Let $\mathcal L$ be the subgroup of $\bfZ_{\mathfrak g}$ generated by
$\lambda_1, \cdots, \lambda_d$. Then $\pi_1(S) \cong \bfZ_{\mathfrak g}/\mathcal L$.
Obviously $\bfZ_{\mathfrak g}/\mathcal L$ is not trivial if $\ell > 1$. Thus we proved the
following.

\begin{proposition} Let $S$ be a compact connected toric Sasaki manifold 
associated with a toric diagram of height $\ell > 1$. Then $S$ is not simply
connected.
\end{proposition}

Note that the converse is not true as the following example shows. Consider the
toric diagram with height $1$ defined by the three normal vectors
$$
 \lambda_1 = \left(\begin{array}{c} 1 \\ 0 \\ 0\end{array}\right), 
\lambda_2 = \left(\begin{array}{c} 1 \\ 2 \\ 1\end{array}\right), 
\lambda_3 = \left(\begin{array}{c} 1 \\ 3 \\ 4 \end{array}\right).
$$
The resulting Sasaki manifold is the Lense space with a different $\bfZ_5$-action
from the above example with $\ell = 5$.
Note also for example $Y^{p,q}$ is not simply connected unless $p$ and $q$ are
relatively prime.

\medskip

\noindent
{\it Proof of Theorem \ref{Main4}} : Put $n = k+3$. 
We
construct diagrams of height $1$ with either $|p_{i+1} - p_i| = 1$ or $|q_{i+1} - q_i| = 1$
such that $(p_1, q_1), \cdots, (p_n, q_n), (p_1,q_1)$ form a convex polytope with $n$ vertices
and that they generate $\bfZ^2$. Then $S$ is simply connected since $\mathcal L = \bfZ^3$.
By another theorem of Lerman \cite{Lerman02} we know that $b_2(S) = n - 3 = k$. 
It follows from the classification of five dimensional simply connected spin manifolds
with $T^3$-action (\cite{BoGaOr06}) that $S = S^5 \sharp k(S^2 \times S^5)$.

There are many ways to construct such examples. For instance
if $k=2r$ so that $n = 2r + 3$ then take
$$
 \left(\begin{array}{c}
p_0 \\ q_0\end{array}\right) 
=  \left(\begin{array}{c}0 \\ 0\end{array}\right),\  
\left(\begin{array}{c}p_1 \\q_1\end{array}\right) 
=  \left(\begin{array}{c}1\\1\end{array}\right), $$

$$\cdots,  
\left(\begin{array}{c}p_r\\ q_r\end{array}\right) =  
\left(\begin{array}{c}r\\ \frac{r(r+1)}2\end{array}\right), 
 \left(\begin{array}{c}p_{r+1}\\q_{r+1}\end{array}\right) 
 =  \left(\begin{array}{c}r+1\\ \frac{(r+1)(r+2)}2 + s\end{array}\right),$$

$$  \left(\begin{array}{c}p_{r+2}\\q_{r+2}\end{array}\right) =  
\left(\begin{array}{c}r\\ \frac{(r+1)(r+2)}2 + s -1\end{array}\right), \cdots,  $$

$$
\left(\begin{array}{c}p_{2r+1}\\ q_{2r+1}\end{array}\right)
=  \left(\begin{array}{c}1\\ \frac{(r+1)(r+2)}2 + s - \frac{r(r+1)}2\end{array}\right), 
\left(\begin{array}{c}p_{2r+2}\\ q_{2r+2}\end{array}\right) = \left(\begin{array}{c}0\\1\end{array}\right).$$
For different values of $s$ they give inequivalent toric diagrams because they have
different areas.
If $k=2r-1$ so that $n = 2r+2$ then take
$$
  \left(\begin{array}{c}p_0\\q_0\end{array}\right) =  
  \left(\begin{array}{c}0\\0\end{array}\right),
  \  \left(\begin{array}{c}p_1\\q_1\end{array}\right) =  
  \left(\begin{array}{c}1\\1\end{array}\right), 
   \  \left(\begin{array}{c}p_2\\q_2\end{array}\right) =  
  \left(\begin{array}{c}2\\ 3\end{array}\right), $$
  
  $$
  \cdots,
  \left(\begin{array}{c}p_{r-1}\\ q_{r-1}\end{array}\right) =  \left(\begin{array}{c}r-1\\ \frac{(r-1)r}2\end{array}\right), $$
  
  $$
  \left(\begin{array}{c}p_r\\ q_r\end{array}\right) =  
  \left(\begin{array}{c}r\\ \frac{r(r+1)}2 + s\end{array}\right), 
  \left(\begin{array}{c}p_{r+1}\\ q_{r+1}\end{array}\right) =  
  \left(\begin{array}{c}0\\ \frac{(r)(r+1)}2 + s +1\end{array}\right),$$
  
 $$ \left(\begin{array}{c}p_{r+3}\\ q_{r+3}\end{array}\right) =  
 \left(\begin{array}{c}-r\\ \frac{r(r+1)}2 + s\end{array}\right),  \left(\begin{array}{c}p_{r+4}\\ q_{r+4}\end{array}\right) =  \left(\begin{array}{c}-(r-1)\\ \frac{(r-1)r}2\end{array}\right),$$
 
 $$
 \cdots,  \left(\begin{array}{c}p_{2r}\\q_{2r}\end{array}\right) = 
 \left(\begin{array}{c}-2\\ 3\end{array}\right)
 \left(\begin{array}{c}p_{2r+1}\\q_{2r+1}\end{array}\right) = 
 \left(\begin{array}{c}-1\\0\end{array}\right).$$
 
 \noindent
 Then again different values of $s$ give inequivalent diagrams. This completes the 
 proof of Theorem \ref{Main4}.

\end{document}